\documentclass[amstex,12pt,amssymb]{article}

\usepackage{mathtext}
\usepackage[cp1251]{inputenc}
\usepackage[T2A]{fontenc}
\usepackage[dvips]{graphicx}
\usepackage{amsmath}
\usepackage{amssymb}
\usepackage{amsxtra}
\usepackage{latexsym}
\usepackage{ifthen}

\newcounter{lemma}[section]

\newcounter{corol}[section]

\newcounter{rem}[section]

\newcounter{theo}[section]

\newcounter{propo}[section]

\numberwithin{equation}{section}

\begin{document}

\markboth{\centerline{E. SEVOST'YANOV}}{\centerline{ON THE LOCAL
BEHAVIOR...}}

\def\cc{\setcounter{equation}{0}
\setcounter{figure}{0}\setcounter{table}{0}}

\overfullrule=0pt


\author{{E. SEVOST'YANOV}\\}

\title{
{\bf ON THE LOCAL BEHAVIOR OF THE MAPPINGS WITH NON--BOUNDED CHARACTERISTICS}}

\date{\today}
\maketitle

\begin{abstract} The present paper is devoted to the study of
space mappings, which are more general than quasiregular mappings.
The questions of the behavior of differentiable mappings having the
so--called $N,$ $N^{-1},$ $ACP$ and $ACP^{-1}$ -- properties are
studied in the work. Under some additional conditions, it is showed
that the modulus of such mappings $f$ can be more than each degree
of logarithmic function at every neighborhood of the isolated
essential singularity of $f.$
\end{abstract}

\bigskip
{\bf 2010 Mathematics Subject Classification: Primary 30C65;
Secondary 30C62}

\section{Introduction}

Here are some definitions. Everywhere below, $D$ is a domain in
${\Bbb R}^n,$ $n\ge 2,$ $m$ be a measure of Lebesgue in ${\Bbb
R}^n,$ and ${\rm dist\,}(A,B)$ is the Euclidean distance between the
sets $A$ and $B$ in ${\Bbb R}^n.$ A mapping $f:D\rightarrow {\Bbb
R}^n$ is said to be a {\it discrete} if the pre-image $f^{-1}(y)$ of
any point $y\,\in\,{\Bbb R}^n$ consists of isolated points, and an
{\it open} if the image of any open set $U\subset D$ is open in
${\Bbb R}^n.$ The notation $f:D\rightarrow {\Bbb R}^n$ assumes that
$f$ is continuous on its domain. In what follows, a mapping $f$ is
supposed to be orientation preserving, i.e. the topological index
$µ(y, f,G)$ is greater than zero for an arbitrary domain $G\subset
D,$ $\overline{G}\subset D$ and an arbitrary $y\in f(G)\setminus
f(\partial G),$ (see, for example, \S\,2 of the Ch. II in
\cite{Re$_2$}). Let $f:D\rightarrow {\Bbb R}^n$ be an arbitrary
mapping and suppose that there is a domain $G\subset D,$
$\overline{G}\subset D,$ for which $
f^{\,-1}\left(f(x)\right)=\left\{x\right\}.$ Then the quantity
$\mu(f(x), f, G),$ which is referred to as the local topological
index, does not depend on the choice of the domain $G$ and is
denoted by $i(x, f).$ In what follows $(x,y)$ denotes the standard
scalar multiplication of the vectors $x,y\in {\Bbb R}^n,$ ${\rm
diam\,}A$ is Euclidean diameter of the set $A\subset {\Bbb R}^n,$
$$B(x_0, r)=\left\{x\in{\Bbb R}^n: |x-x_0|< r\right\}\,,\quad {\Bbb B}^n
:= B(0, 1)\,,$$ $$S(x_0,r) = \{ x\,\in\,{\Bbb R}^n :
|x-x_0|=r\}\,,\quad{\Bbb S}^{n-1}:=S(0, 1)\,,$$ $\omega_{n-1}$
denotes the quare of the unit sphere ${\Bbb S}^{n-1}$ in ${\Bbb
R}^n,$ $\Omega_{n}$ is a volume of the unit ball ${\Bbb B}^{n}$ in
${\Bbb R}^n.$ Given a mapping $f:D\rightarrow {\Bbb R}^n,$ a set
$E\subset D,$ and a point $y\in {\Bbb R}^n$ we define the
multiplicity function $N(y, f, E)$ as the number of pre-images of
$y$ in $E,$ that is, $$N(y, f, E) = {\rm card}\in \left\{x \in E:
f(x) = y\right\}\,.$$

Recall that a mapping $f:D\rightarrow {\Bbb R}^n$ is said to have
the {\it $N$ -- property (of Luzin)} if
$m\left(f\left(S\right)\right)=0$ whenever $m(S)=0$ for all such
sets $S\subset{\Bbb R}^n.$ Similarly, $f$ has the {\it $N^{-1}$ --
property} if $m\left(f^{\,-1}(S)\right)=0$ whenever $m(S)=0.$

We write $f\in W^{1,n}_{loc}(D),$ iff all of the coordinate
functions $f_j,$ $f=(f_1,\ldots,f_n),$ have the partitional
derivatives which are locally integrable in the degree $n$ in $D.$

Recall that a mapping $f:D\rightarrow {\Bbb R}^n$ is said to be {\it
a mapping with bounded distortion}, if the following conditions
hold:

\noindent 1) $f\in W_{loc}^{1,n},$

\noindent 2) a Jacobian $J(x,f):={\rm det\,}f^{\,\prime}(x) $ of the
mapping $f$ at the point $x\in D$ preserves the sign almost
everywhere in $D,$

\noindent 3) $\Vert f^{\,\prime}(x) \Vert^n \le K \cdot |J(x,f)|$ at
a.e. $x\in D$ and some constant $K<\infty,$ where $$\Vert
f^{\,\prime}(x)\Vert:=\sup\limits_{h\in {\Bbb R}^n:
|h|=1}|f^{\,\prime}(x)h|\,,$$
see., e.g., $\S\, 3$ Ch. I in \cite{Re$_2$}, or definition 2.1 of
the section 2 Ch. I in \cite{Ri}.

Active investigations of the mappings with bounded distortion were
started by Yu.G. Reshetnyak. In particular, he has proved that the
mappings $f$ with bounded distortion are open and discrete, see
Theorems 6.3 and 6.4, $\S\, 6,$ Ch. II in \cite{Re$_2$}, are
differentiable a.e., see Theorem 4 in \cite{Re$_1$}, in and have $N$
-- property, see Theorem 6.2 Ch. II in \cite{Re$_2$}. From other
hand, the $N^{-1}$ -- property of the mappings with bounded
distortion was proved by B. Bojarski and T. Iwaniec, see Theorem 8.1
in \cite{BI}.

We recall that an isolated point $x_0$ of the boundary $\partial D$
of a domain $D$ in ${\Bbb R}^n$  is said to be a {\it removable
singularity} if there is a finite limit $\lim\limits_{x\rightarrow
x_0}\,f(x).$ If $f(x)\rightarrow \infty$ as $x\rightarrow x_0,$ then
$x_0$ is referred to as a {\it pole}. An isolated point $x_0$ of
$\partial D$ is called an {\it essential singularity} of a mapping
$f:D\rightarrow {\Bbb R}^n$ if the limit $\lim\limits_{x\rightarrow
x_0}\,f(x)$ does not exist.

In 1972, in the work of J.~V\"{a}is\"{a}l\"{a} was proved the
following, see e.g. Theorem 4.2 in \cite{Va$_2$}.

\medskip
{\bf Statement 1.} {\sl Let $b\in D$ and
$f:D\setminus\{b\}\rightarrow {\Bbb R}^n$ be a mapping with bounded
distortion. Suppose that there exists a number $\delta>0$ such that
\begin{equation}\label{eq1}
|f(x)|\le C |x-b|^{-p}\,,
\end{equation}
at every $x\in B(b, \delta)\setminus\{b\}$ and some positive
constants $p>0$ and $C>0.$ Then $b$ is a removable singularity or a
pole of the mapping $f.$}

\medskip
A goal of the present paper is a proof of the analogue of the
statement 1 for more general classes of mappings of finite length
distortion, including the classes of mappings with bounded
distortion. Mappings with finite length distortion were introduced
by O. Martio, V. Ryazanov, U. Srebro and E. Yakubov in 2002, see
e.g. in the work \cite{MRSY$_1$}, or Chapter 8 in \cite{MRSY$_2$}.
The considering of it is actually in the connection with the study
of the so--called mappings with finite distortion, which are
actively investigated at the last time, see e.g. Chapter 20 in
\cite{AIM} or Chapter 6 in \cite{IM}. In this connection, see also
the works \cite{BGMV}, \cite{Mikl}, \cite{Sal} and \cite{UV}.

A curve $\gamma$ in ${\Bbb R}^n$ is a continuous mapping $\gamma
:\Delta\rightarrow{\Bbb R}^n$ where $\Delta$ is an interval in
${\Bbb R} .$ Its locus $\gamma(\Delta)$ is denoted by $|\gamma|.$
Given a family of curves $\Gamma$ in ${\Bbb R}^n ,$ a Borel function
$\rho:{\Bbb R}^n \rightarrow [0,\infty]$ is called {\it admissible}
for $\Gamma ,$ abbr. $\rho \in {\rm adm}\, \Gamma ,$ if curvilinear
integral of the first type $\int\limits_{\gamma} \rho(x)|dx|$
satisfies the condition $$\int\limits_{\gamma} \rho(x)|dx| \ge 1$$
for each $\gamma\in\Gamma.$ The {\it modulus} $M(\Gamma )$ of
$\Gamma$ is defined as
$$M(\Gamma) =\inf\limits_{ \rho \in {\rm adm}\, \Gamma}
\int\limits_{{\Bbb R}^n} \rho^n(x) dm(x)$$ interpreted as $+\infty$
if ${\rm adm}\, \Gamma = \varnothing .$ The properties of the above
modulus are analogous to the properties of the measure of Lebesgue
$m$ in ${\Bbb R}^n.$ Namely, a modulus of the empty family equals to
zero, $M(\varnothing)=0,$ a modulus has a property of monotonicity
by the relation to families of curves $\Gamma_1$ and $\Gamma_2:$
$\Gamma_1\subset\Gamma_2\Rightarrow M(\Gamma_1)\le M(\Gamma_2),$ and
has a property of subadditivity,
$M\left(\bigcup\limits_{i=1}^{\infty}\Gamma_i\right)\le
\sum\limits_{i=1}^{\infty}M(\Gamma_i),$ see Theorem 6.2 in
\cite{Va$_1$}.

We say that a property $P$ holds for {\it almost every (a.e.)}
curves $\gamma$ in a family $\Gamma$ if the subfamily of all curves
in $\Gamma $ for which $P$ fails has modulus zero.

If $\gamma :\Delta\rightarrow{\Bbb R}^n$ is a locally rectifiable
curve, then there is the unique increasing length function
$l_{\gamma}$ of $\Delta$ onto a length interval $\Delta
_{\gamma}\subset{\Bbb R}$ with a prescribed normalization $l
_{\gamma}(t_0)=0\in\Delta _{\gamma},$ $t_0\in\Delta,$ such that $l
_{\gamma}(t)$ is equal to the length of the subcurve $\gamma
|_{[t_0,t]}$ of $\gamma$ if $t>t_0,$ $t\in\Delta ,$ and $l
_{\gamma}(t)$ is equal to $-l(\gamma |_{[t,t_0]})$ if $t<t_0,$
$t\in\Delta .$ Let $g: |\gamma |\rightarrow{\Bbb R}^n$ be a
continuous mapping,  and suppose that the curve $\widetilde{\gamma}
=g\circ\gamma$ is also locally rectifiable. Then there is a unique
increasing function $L_{\gamma ,g}: \Delta
_{\gamma}\rightarrow\Delta _{\widetilde{\gamma}}$ such that
$L_{\gamma ,g}\left(l_{\gamma}(t)\right) = l_{\widetilde{\gamma}}(t)
\quad\forall\quad t\in\Delta.$ A curve $\gamma$ in $D$ is called
here a {\it lifting} of a curve $\widetilde{\gamma}$ in ${\Bbb R}^n$
under $f:D\rightarrow {\Bbb R}^n$ if $\widetilde{\gamma} =
f\circ\gamma.$ Recall that $f\in ACP$ if and only if a curve
$\widetilde{\gamma}=f\circ\gamma$ is locally rectifiable for a.e.
curves $\gamma$ in $D,$ and $L_{\gamma , f} $ is absolutely
continuous on closed subintervals of $\Delta_{\gamma}$ for a.e.
curves $\gamma$ in $D.$  We say that a discrete mapping $f$ is {\it
absolute continuous on curves in the inverse direction,} abbr.
$ACP^{-1},$ if for a.e. curves $\widetilde{\gamma}$ a lifting
$\gamma$ of $\widetilde{\gamma},$ $\widetilde{\gamma}=f\circ\gamma,$
is locally rectifiable, and $L^{-1}_{\gamma , f}$ is absolutely
continuous on closed subintervals of $\Delta_{\widetilde{\gamma}}$
for a.e. curves $\widetilde{\gamma}$ in $f(D)$ and for each lifting
$\gamma$ of $\widetilde{\gamma}.$ A mapping $f:D\rightarrow{\Bbb
R}^n$ is said to be of {\it finite length distortion}, abbr. $f\in
FLD$, if $f$ is differentiable a.e. in $D,$ has $N$ -- and $N^{-1}$
-- properties, and $f\in ACP\cap ACP^{-1}.$

\begin{rem}\label{rem1} The notion of the mappings with finite length distortion
can be given in more general case, when $f$ does not supposed to be
a discrete, see e.g. in \cite{MRSY$_1$}, see also section 8.1 in
\cite{MRSY$_2$}. Of course, the above definition is equivalent to
the correspondent general case, see, for instance, section 8.1 and
corollary 8.1 in \cite{MRSY$_2$}, or Corollary 3.14 in
\cite{MRSY$_1$}. In this connection, the word $"$discrete$"$ will be
present in the text, if it is necessary.

For the classes $W_{loc}^{1,n},$ and, in particular, for the
mappings with bounded distortion, the $ACP$ property is well--known
as B. Fuglede's lemma, see, for instance, Theorem 28.2 in
\cite{Va$_1$}. Besides that, the $ACP^{-1}$ property was proved by
E.A. Poletskii for it, see e.g. Lemma 6 in \cite{Pol}. Taking into
account all of the comments given above, we conclude that every
mapping with bounded distortion is a mapping with finite length
distortion, see also Theorem 4.7 in \cite{MRSY$_1$}, or Theorem 8.2
in \cite{MRSY$_2$} in this connection.
\end{rem}

We say that a function ${\varphi}:D\rightarrow{\Bbb R}$ has a {\it
finite mean oscillation} at the point $x_0\in D$, write $\varphi\in
FMO(x_0),$ if
$${\limsup\limits_{\varepsilon\rightarrow 0}}\
\frac{1}{\Omega_n\cdot\varepsilon^n} \int\limits_{B( x_0,
\varepsilon)} |{\varphi}(x)-\overline{{\varphi}}_{\varepsilon}|\
dm(x)<\infty\,,$$
where $\overline{{\varphi}}_{\varepsilon}=
\frac{1}{\Omega_n\cdot\varepsilon^n}\int\limits_{B(
x_0,\,\varepsilon)} {\varphi}(x)\ dm(x).$
Functions of finite mean oscillation were introduced by A. Ignat'ev
and V. Ryazanov in the work \cite{IR}, see also section 11.2 in
\cite{MRSY$_2$}. There are the generalization and localization of
the space $BMO,$ that is bounded mean oscillation functions by F.
John and L. Nirenberg, see for instance \cite{JN}.

\medskip
Set $l\left(f^{\,\prime}(x)\right):=\inf\limits_{h\in {\Bbb R}^n:
|h|=1}|f^{\,\prime}(x)h|.$ Recall that {\it inner dilatation} of the
mapping $f$ at a point $x$ is defined as
$$K_I(x,f)\quad =\quad= \left\{
\begin{array}{rr}
\frac{|J(x,f)|}{{l\left(f^{\,\prime}(x)\right)}^n}, & J(x,f)\ne 0,\\
1,  &  f^{\,\prime}(x)=0, \\
\infty, & {\rm otherwise}
\end{array}
\right.\,.$$ {\it Outher dilatation} of the mapping $f$ at the point
$x$ can be defined as
$$K_O(x,f)\quad =\quad= \left\{
\begin{array}{rr}
\frac{\Vert f^\prime(x)\Vert^n}{|J(x,f)|}, & J(x,f)\ne 0,\\
1,  &  f^{\,\prime}(x)=0, \\
\infty, & {\rm otherwise}
\end{array}
\right.\,.$$ It is well--known that $K_I(x,f)\le K_O^{n-1}(x,f)$
everywhere at the points, where there are well--defined, see for
instance formulae (2.7) and (2.8) of the section 2.1 of Ch. I in
\cite{Re$_2$}. In particular, for the mappings with bounded
distortion we have $K_I(x, f)\le K^{n-1}$ at a.e. $x,$ that follows
from it's definition. The main result of the paper is the following.

\medskip
{\bf Statement $1^{\prime}.$} {\sl Let $b\in D$ and
$f:D\setminus\{b\}\rightarrow {\Bbb R}^n$ be an open and discrete
mapping with finite length distortion. Suppose that there exists
$\delta>0$ such that
\begin{equation}\label{eq2}
|f(x)|\le C \left(\log\frac{1}{|x-b|}\right)^{p}\,,
\end{equation}
at every $x\in B(b, \delta)\setminus\{b\}$ and some constants $p>0$
and $C>0.$ Let there exists a function $Q:D\rightarrow[1, \infty],$
such that $K_I(x,f)\le Q(x)$ a.e. $x\in D$ and $Q(x)\in FMO(b).$
Then a point $b$ is a removable singularity, or a pole of the
mapping $f.$}

\begin{rem}\label{rem2}
Note that the condition (\ref{eq2}) is stronger than the requirement
(\ref{eq1}), and that from the Statement $1^{\prime}$ it follows the
Statement 1.
\end{rem}

\section{Preliminaries. The main Lemma}
\setcounter{equation}{0}

The following definitions can be found in the monograph \cite{Ri},
Ch. II, Section. 3, see also section  3.11 in \cite{Va$_2$}. Let
$f:D \rightarrow {\Bbb R}^n$ be an arbitrary mapping,
$\beta:[a,\,b)\rightarrow {\Bbb R}^n$ is a path and
$x\in\,f^{\,-1}\left(\beta(a)\right).$ A path $\alpha:
[a,\,c)\rightarrow D$ is called a {\it maximal $f$ -- lifting} of
$\beta$ starting at $x,$ if $(1)\quad\alpha(a)=x;$ $(2)\quad
f\circ\alpha=\beta|_{[a,\,c)};$ $(3)$\quad if $c<c^{\prime}\,\le b,$
then there does not exist a path $\alpha^{\prime}: [a,
c^{\prime})\rightarrow D$ such that $\alpha =
\alpha^{\prime}|_{[a,\,c)}$ and $f\circ
\alpha^{\,\prime}=\beta|_{[a,\,c^{\prime})}.$

Let $x_1,\ldots,x_k$ be $k$ different points of
$f^{-1}\left(\beta(a)\right)$ and let
$$m = \sum\limits_{i=1}^k i(x_i,\,f).$$
We say that the sequence $\alpha_1,\dots,\alpha_m$ is a {\it maximal
sequence of $f$ -- lifting of $\beta$ starting at points
$x_1,\ldots,x_k,$} if

$(a)$\quad each $\alpha_j$ is a maximal $f$ -- lifting of $\beta,$

$(b)\quad {\rm card}\,\left\{j:a_j(a)=x_i\right\}= i(x_i,\,f),\quad
1\le i\le k\,,$

$(c)\quad {\rm card}\,\left\{j:a_j(t)=x\right\}\le i(x,\,f)$ for all
$x\in D$ and for all $t.$

Let $f$ be a discrete open mapping and $x_1,\ldots,x_k$ are distinct
points in $f^{\,-1}\left(\beta(a)\right).$ Then $\beta$ has a
maximal sequence of $f$ -- liftings starting at points
$x_1,\ldots,x_k,$ see Theorem 3.2 Ch. II in \cite{Ri}. The following
statement was proved by author, see for instance Theorem 3.1 in
\cite{Sev$_3$}.

\begin{propo}{}\label{pr1}{\sl\,
Let $f:D\rightarrow {\Bbb R}^n$ be a discrete and an open mapping
with finite length distortion, $\Gamma$ a path family in $D,$
$\Gamma^{\,\prime}$ a path family in ${\Bbb R}^n$ and $m$ a positive
integer such that the following is true. Suppose that for every path
$\beta$ in $\Gamma^{\,\prime}$ there are paths
$\alpha_1,\ldots,\alpha_m$ in $\Gamma$ such that $f\circ
\alpha_j\subset \beta$ for all $j=1,\ldots,m$ and such that for
every $x\in D$ and all $t$ the equality $\alpha_j(t)=x$ holds for at
most $i(x,f)$ indices $j.$ Then
\begin{equation}\label{eq5}
M(\Gamma^{\,\prime} )\quad\le\quad \frac{1}{m}\quad\int\limits_D
K_I\left(x,\,f\right)\cdot \rho^n (x)\ \ dm(x)
\end{equation}
for every $\rho \in {\rm \,adm}\,\Gamma.$}
\end{propo}

In particular, the Proposition \ref{pr1} generalizes the
corresponding result of O. Martio, V. Ryazanov, U. Srebro and E.
Yakubov, that is the above statement under $m=1,$ see for instance
Theorem 6.10 in \cite{MRSY$_1$}, of Theorem 8.6 in \cite{MRSY$_2$}.
By $\Gamma(E,F,D)$ we denote the family of all curves
$\gamma:[a,b]\rightarrow{\Bbb R}^n$ connecting $E$ and $F$ in $D,$
i.e. $\gamma(a)\in E,$ $\gamma(b)\in F$ and $\gamma(t)\in D$ as
$t\in (a, b).$ A compact set $G\subset {\Bbb R}^n$ is said to be
{\it a set of capacity zero,} write ${\rm cap\,}G =0,$ if there
exists $T\subset {\Bbb R}^n,$ such that $M(\Gamma(T, G, {\Bbb
R}^n))=0,$ see, for instance, Section 2 of Ch. III and Proposition
10.2 of Ch. II in \cite{Ri}. By definition, an arbitrary set has a
zero capacity if and only if every it's compact subset has a zero
capacity. The sets of capacity zero are totally disconnected, i.e.,
the condition ${\rm cap\,}G =0$ implies that ${\rm Int\,}G=0,$ see
e.g. Corollary 2.5 of Ch. III in \cite{Ri}. Open set $U\subset D,$
$\overline{U}\subset D,$ is said to be a {\it normal neighborhood}
of the point $x\in D$ under the mapping $f:D\rightarrow {\Bbb R}^n,$
iff $U\cap f^{\,-1}\left(f(x)\right)=\left\{x\right\}$ and $\partial
f(U)=f(\partial U),$ see e.g. Section 4 of Ch. I in \cite{Ri}.

\begin{propo}\label{pr2}{\sl\,
Let $f:D\rightarrow {\Bbb R}^n$ be an open discrete mapping, then
for every $x\in D$ there exists $s_x$ such that, for every $s\in (0,
s_x),$ the $x$ -- component of $f^{-1}\left(B(f(x), s)\right),$
denoted by $U(x,f,s),$ is a normal neighborhood of $x$ under $f,$
$f\left(U(x,f,s)\right)=B(f(x), s)$ and ${\rm
diam\,}U(x,f,s)\rightarrow 0$ as $s\rightarrow 0,$ see, for
instance, Lemma 4.9 of Ch. I in \cite{Ri}.}
\end{propo}

The main tool under the proof of the basic results of the present
work is the following

\begin{lemma}{}\label{lem1}{\sl\, Let $b\in D$ and
$f:D\setminus\{b\}\rightarrow {\Bbb R}^n$ be an open and a discrete
mapping with finite length distortion. Suppose that there exists
$\delta>0$ such that
\begin{equation}\label{eq3}
|f(x)|\le C \left(\log\frac{1}{|x-b|}\right)^{p}\,,
\end{equation}
at every $x\in B(b, \delta)\setminus\{b\}$ and some constants $p>0$
and $C>0.$ Follow, suppose that there exist a measurable function
$Q:D\rightarrow [1, \infty],$ numbers $\varepsilon_0>0,$
$\varepsilon_0<{\rm dist\,}\left(b,
\partial D\setminus\{b\}\right),$ $A>0$ and a Borel
function $\psi(t):[0, \varepsilon_0]\rightarrow
(0, \infty)$ such that $K_I(x,f)\le Q(x)$ a.e., and
\begin{equation}
\label{eq4}
\int\limits_{\varepsilon<|x-b|<\varepsilon_0}Q(x)\cdot\psi^n(|x-x_0|)
\ dm(x)\le \frac{A\cdot I^n(\varepsilon,
\varepsilon_0)}{\left(\log\log\frac{1}{\varepsilon}\right)^{n-1}}
\qquad \forall\quad\varepsilon\in(0,\varepsilon_0/2)\,,
\end{equation} где
\begin{equation} \label{eq11}
0<I(\varepsilon, \varepsilon_0)
=\int\limits_{\varepsilon}^{\varepsilon_0}\psi(t)dt < \infty \qquad
\forall\quad\varepsilon \in(0, \varepsilon_0)\,.\end{equation} Then
a point $b$ is a removable singularity, or a pole of the mapping
$f.$}
\end{lemma}

{\it Proof.\,} Suppose the contrary, i.e., a point $b$ is an
essential singularity of $f.$ Without loss of generality, we can
consider that $b=0$ and $C=1.$ In this case, there exists $R>0,$
such that
\begin{equation}\label{eq10}
f\left(S(0, \delta)\right)\subset B(0, R)\,.
\end{equation}
Since $b=0$ is an essential singularity of $f,$ from the conditions
(\ref{eq5}), (\ref{eq4}) and another author's result, see Lemma 3.1,
Lemma 5.1 and  Theorem 6.5 in \cite{Sev$_4$}, we have $$N\left(y, f,
B(0, \delta)\setminus \{0\}\right)=\infty$$ for all $y\in {\Bbb
R}^n\setminus E,$ where ${\rm cap\,}E=0.$ Since $E$ of zero
capacity, ${\Bbb R}^n\setminus E$ is unbounded. Thus, there exists
$y_0\in {\Bbb R}^n\setminus \left(E\cup B(0, R)\right).$

Let $k_0>\frac{4Ap^{n-1}}{\omega_{n-1}},$ $k_0\in {\Bbb N}.$ Since
$N\left(y_0, f, B(0, \delta)\setminus \{0\}\right)=\infty,$ there
exist the points $x_1,\ldots,x_{k_0}\in f^{-1}(y_0),$
$x_1,\ldots,x_{k_0}\in B(0, \delta)\setminus\{0\}.$ By Proposition
\ref{pr2}, there exists $r>0,$ such that every point $x_j,$
$j=1,\ldots,k_0,$ has a normal neighborhood $U_j:=U(x_j, f , r)$
with $\overline{U_l}\cap\overline{U_m}=\varnothing$ at all $l\ne m,
$ $l, m\in {\Bbb N},$ $1\le l\le k_0$ and $1\le m\le k_0.$

Set $d:=\min\left\{\varepsilon_0, {\rm dist\,}\left(0,
\overline{U_1}\cup\ldots\cup \overline{U_{k_0}}\right)\right\}.$ Let
$a\in (0, d)$ and  $V:=B(0, \delta)\setminus\overline{B(0, a)}.$ By
(\ref{eq3}), taking into account that $\partial f(V)\subset
f(\partial V)$ and $C=1,$ we have
\begin{equation}\label{eq6}
f(V)\subset B\left(0, \left(\log\frac{1}{a}\right)^p\right)\,.
\end{equation}
Since $z_0:=y_0+re\in \overline{B(y_0, r)}=f\left(\overline{U(x_j,
f, r)}\right),$ $j=1,\ldots, k_0,$ we have $z_0\in f(V).$ Thus there
exists a sequence of the points
$\widetilde{x_1},\ldots,\widetilde{x_{k_0}},$ $\widetilde{x_j}\in
\overline{U_j},$ $1\le j\le k_0,$ such that
$f(\widetilde{x_j})=z_0.$ Note that $k_0\le\sum\limits_{j=1}^{k_0}
i(\widetilde{x_j},\,f)=m^{\,\prime}.$ Let $H$ be a hemisphere
$H=\left\{e\in {\Bbb S}^{n-1}: (e, y_0)>0\right\},$
$\Gamma^{\,\prime}$ be a curve's family $\beta:\left[r,
\left(\log\frac{1}{a}\right)^p\right)\rightarrow {\Bbb R}^n$ of the
type $\beta(t)=y_0+te,$ $e\in H,$ and $\Gamma$ be a sequence of
maximal $f$ -- liftings of $\beta$ under the mapping $f$ in $V,$
starting at the points $\widetilde{x_1},\ldots,\widetilde{x_{k_0}},$
$\widetilde{x_j}\in \overline{U_j},$ $1\le j\le k_0,$ consisting
from $m^{\,\prime}$ curves, where $m^{\,\prime}
=\sum\limits_{j=1}^{k_0} i(\widetilde{x_j},\,f).$ Such a sequence
exists by Theorem 3.2 of Ch. II in \cite{Ri}.  By Proposition
\ref{pr1},
\begin{equation}\label{eq9} M(\Gamma^{\,\prime} )\le
\frac{1}{m^{\,\prime}}\quad\int\limits_D K_I\left(x,\,f\right)\cdot
\rho^n (x) dm(x)\le \frac{1}{k_0}\int\limits_D
K_I\left(x,\,f\right)\cdot \rho^n (x) dm(x)
\end{equation}
for every $\rho \in {\rm }\,{\rm adm}\,\Gamma.$

Given $e\in H,$ we show that, for every curve $\beta=y_0+te$ and
it's maximal lifting $\alpha(t):[r, c)\rightarrow V$ starting at the
point $\widetilde{x_{j_0}},$ $\alpha\in \Gamma,$ $1\le j_0\le k_0,$
there exists a sequence  $r_k\in [r, c)$ with $r_k\rightarrow c-0$
as $k\rightarrow \infty$ such that ${\rm dist\,}(\alpha(r_k),
\partial V)\rightarrow 0$ as $k\rightarrow \infty.$
Suppose the contrary, i.e., there exists $e_0\in H,$ such that
$\alpha(t),$ $t\in [r, c),$ is a maximal lifting of
$\beta=y_0+te_0,$ and $\alpha(t)$ lies inside of $V$ with it's
closure. Let $C(c,\,\alpha(t))$ denotes a cluster set of $\alpha$ as
$t\rightarrow c-0.$ For every $x\in C(c,\,\alpha(t))$ there exists a
sequence $t_k\rightarrow \infty$ such that
$x=\lim\limits_{k\rightarrow \infty}\alpha(t_k).$ Since $f$ is
continuous and $C(c,\,\alpha(t))\subset V$ by the assumption, we
have $f(x)=f(\lim\limits_{k\rightarrow \infty}\alpha(t_k))=
\lim\limits_{k\rightarrow \infty}\beta(t_k)= \beta(c),$ from what it
follows that $f$ is a constant on $C(c,\,\alpha(t)).$ Since $f$ is a
discrete and a set $C(c,\,\alpha(t))$ is connected, we have
$C(c,\,\alpha(t))=p_1\in V.$ Let $c\ne
b_0:=\left(\log\frac{1}{a}\right)^p.$ In this case, we can construct
a lifting  $\alpha^{\,\prime}$ of $\beta$ started at $p_1.$ Uniting
the liftings $\alpha$ and $\alpha^{\,\prime},$ we obtain another
maximal lifting $\alpha^{\,\prime\prime}$ of $\beta$ starting at the
point $\widetilde{x_{j_0}},$ that contradicts to the maximal
property of the first lifting $\alpha.$ Thus, $c=b_0$ and hence
$\alpha$ can be extend to closed curve defined on the segment
$\left[r, \left(\log\frac{1}{a}\right)^p\right]$ (we don't change
the notion of the extended curve). Then, at every $t\in \left[r,
\left(\log\frac{1}{a}\right)^p\right],$ we have
$\beta(t)=f(\alpha(t))\subset f(V).$ In particular, by (\ref{eq6})
\begin{equation}\label{eq7}
z_1:=y_0+\left(\log\frac{1}{a}\right)^pe_0\in f(V)\subset B\left(0,
\left(\log\frac{1}{a}\right)^p\right)\,.
\end{equation}
However, since $e_0\in H,$ we have
$$|z_1|=\left|y_0+\left(\log\frac{1}{a}\right)^pe_0\right|=\sqrt{|y_0|^2 +
2\left(y_0,
\left(\log\frac{1}{a}\right)^pe_0\right)+\left(\log\frac{1}{a}\right)^{2p}}\ge$$
\begin{equation}\label{eq8}
\ge \sqrt{|y_0|^2 + \left(\log\frac{1}{a}\right)^{2p}}\ge
\left(\log\frac{1}{a}\right)^{p}\,.
\end{equation}
However, the relation (\ref{eq8}) contradicts to (\ref{eq7}), which
disproves the assumption that $\alpha(t)$ consists in the set $V$
with it's closure.  Consequently, ${\rm dist\,}(\alpha(r_k),
\partial V)\rightarrow 0$ as $k\rightarrow c-0$ and some sequence $r_k\in [r, c)$
such that $r_k\rightarrow c-0$ as $k\rightarrow \infty.$

Note that the situation when ${\rm dist\,}(\alpha(r_k), S(0,
\delta))\rightarrow 0$ as $k\rightarrow c-0$ is excluded. In fact,
suppose that there exist $p_2\in S(0, \delta)$ and a sequence $k_l,$
$l\in {\Bbb N},$ such that $\alpha(r_{k_l})\rightarrow p_2$ as
$l\rightarrow \infty.$ By the continuously of $f$ we have that
$\beta(r_{k_l})\rightarrow f(p_2)$ as $l\rightarrow \infty,$ that is
impossible by (\ref{eq10}), because for every  $e\in H$ and $t\in
\left[r, \left(\log\frac{1}{a}\right)^p\right)$ we have
$|\beta(t)|=|y_0+te|=\sqrt{|y_0|^2 + 2t(y_0, e)+t^2}\ge |y_0|> R$ by
the choosing of $y_0.$

It follows from above that there exists a sequence $r_k\in [r, c)$
such that $r_k\rightarrow c-0$ as $k\rightarrow \infty$ and
$\alpha(r_k)\rightarrow p_3\in S(0, a).$ Besides that, every such a
curve $\alpha\in \Gamma$ intersects the sphere $S(0, d)$ because
$\alpha$ has a start outside of $B(0, d).$ Consider the function
$$\rho_{a}(x)= \left\{
\begin{array}{rr}
\psi(|x|)/I(a, d), &   x\in B(0,d)\setminus B(0,a),\\
0,  &  x\in {\Bbb R}^n \setminus  \left(B(0,d)\setminus
B(0,a)\right)
\end{array}
\right.\,, $$ where $I(a, d)$ is defined as in (\ref{eq11}) and
$\psi$ be a function from the condition of Lemma. Since $\psi(t)>0,$
$I(a, d)>0$ for every $0<a<d.$ Thus, a function $\rho_{a}(x)$ which
is given above is well--defined. Note that a function $\rho_{a}(x)$
is Borel, moreover, since $\rho_{a}(x)$ is a radial function, by the
above properties of curves of $\Gamma$ and by Theorem 5.7 in
\cite{Va$_1$}, for every curve $\alpha\in \Gamma$ we have
$$\int\limits_{\alpha}\rho_a(x)|dx|\ge\frac{1}{I(a, d)}\int\limits_a^d \psi(t)dt= 1\,,$$
i.e., $\rho_{a}(x)\in {\rm adm\,}\Gamma.$ Thus, from (\ref{eq4}) and
(\ref{eq9}) we have
$$M(\Gamma^{\,\prime} )\quad\le\quad \frac{1}{k_0\cdot I^n(a,
d)}\quad\int\limits_{a<|x|<d} K_I\left(x,\,f\right)\cdot \psi^n
(|x|)dm(x)\le$$ $$\le \frac{I^n(a, \varepsilon_0)}{k_0\cdot I^n(a,
d)\cdot I^n(a, \varepsilon_0)}\quad\int\limits_{a<|x|<\varepsilon_0}
Q(x)\cdot \psi^n (|x|)dm(x)=$$$$=\left(1+\frac{I(d,
\varepsilon_0)}{I(a, d)}\right)^n\frac{1}{k_0\cdot I^n(a,
\varepsilon_0)}\quad\int\limits_{a<|x|<\varepsilon_0} Q(x)\cdot
\psi^n (|x|)dm(x)\le$$ $$\le \frac{2}{k_0\cdot I^n(a,
\varepsilon_0)}\quad\int\limits_{a<|x|<\varepsilon_0} Q(x)\cdot
\psi^n (|x|)dm(x)$$ at every $a\in (0, d_1)$ and some $d_1,$ $d_1\le
d,$ because by (\ref{eq4}), $I^n(a,d)\rightarrow\infty$ as
$a\rightarrow \infty.$ Now, from (\ref{eq4}), we conclude that
\begin{equation}\label{eq12}
M(\Gamma^{\,\prime} )\le
\frac{2A}{k_0\left(\log\log\frac{1}{a}\right)^{n-1}}
\end{equation}
for $a\in (0, d_1).$ From other hand, by section 7.7 in
\cite{Va$_1$},
\begin{equation}\label{eq13}
M(\Gamma^{\prime})=\frac{1}{2}\frac{\omega_{n-1}}
{\left(\log\frac{\left(\log\frac{1}{a}\right)^p}{r}\right)^{n-1}}\,.
\end{equation}
By (\ref{eq12}) and (\ref{eq13}) we have
$$\frac{1}{2}\frac{\omega_{n-1}}
{\left(\log\frac{\left(\log\frac{1}{a}\right)^p}{r}\right)^{n-1}}\le
\frac{2A}{k_0\left(\log\log\frac{1}{a}\right)^{n-1}},$$ from what
$$\left(\log\left(\frac{\left(\log\frac{1}{a}\right)^p}{r}\right)^{
\left(\frac{2}{\omega_{n-1}}\right)^{\frac{1}{n-1}}}\right)^{n-1}\ge
\left(\log\left(\log\frac{1}{a}\right)^{\left(\frac{k_0}{2A}\right)^{\frac{1}{n-1}}}\right)^{n-1}\,,
$$
$$\left(\frac{\left(\log\frac{1}{a}\right)^p}{r}\right)^{
\left(\frac{2}{\omega_{n-1}}\right)^{\frac{1}{n-1}}}\ge
\left(\log\frac{1}{a}\right)^{\left(\frac{k_0}{2A}\right)^{\frac{1}{n-1}}}\,,
$$
$$\frac{1}{r^{{
\left(\frac{2}{\omega_{n-1}}\right)^{\frac{1}{n-1}}}}}\ge
\left(\log\frac{1}{a}\right)^{{\left(\frac{k_0}{2A}\right)^{\frac{1}{n-1}}}-p{
\left(\frac{2}{\omega_{n-1}}\right)^{\frac{1}{n-1}}}}\,.$$ Since by
the choosing $k_0>\frac{4Ap^{n-1}}{\omega_{n-1}},$ in the
right--hand part of the above relation the logarithmic function
presents in some positive degree. Letting to the limit as
$a\rightarrow 0$ in the right--hand part of it we obtain that
$$\frac{1}{r^{{
\left(\frac{2}{\omega_{n-1}}\right)^{\frac{1}{n-1}}}}}\ge
\infty\,,$$ that is impossible. The contradiction obtained above
disproves the assumption that $b=0$ is an essential singularity of
$f.$ $\Box$

\medskip
The next statement follows directly from Lemma 5 in \cite{Sev$_1$}
as $\psi(t)=\frac{1}{t\log\frac{1}{t}}$ and from the estimate
(\ref{eq5}) at $m=1.$

\begin{propo}{}\label{pr3}
{\sl\, Let $b\in D$ and $f:D\rightarrow {\Bbb R}^n$ be an open and a
discrete mapping with finite length distortion. Suppose that there
exist a measurable function $Q:D\rightarrow [1, \infty],$ the
numbers $\varepsilon_0>0,$ $\varepsilon_0<{\rm dist\,}\left(b,
\partial D\right),$ and $A>0$ such that $K_I(x,f)\le Q(x)$ a.e., such that the
relations (\ref{eq4}) and (\ref{eq11}) hold as
$\psi(t)=\frac{1}{t\log\frac{1}{t}},$ i.e.,
\begin{equation}\label{eq14}
\int\limits_{\varepsilon<|x-b|<\varepsilon_0}\frac{Q(x)}{|x-b|^n
\log^n\frac{1}{|x-b|}}\ dm(x)\le A\cdot\log{\frac{\log{\frac{1}
{\varepsilon}}}{\log{\frac{1}{\varepsilon_0}}}} \qquad
\forall\quad\varepsilon\in(0,\varepsilon_0)\,.
\end{equation}
Then for every $x\in B(b, \varepsilon_0)$
\begin{equation}\label{eq15}
|f(x)-f(b)|\le
\frac{\alpha_n(1+R^2)}{\delta}\left\{\frac{\log\frac{1}{\varepsilon_0}}
{\log\frac{1}{|x-b|}}\right\}^{\beta_n}\,,
\end{equation}
where $\alpha_n$ and
$\beta_n=\left(\frac{\omega_{n-1}}{A}\right)^{1/(n-1)}$ depend only
on $n,$ and $\delta$ depends only on $R.$ }
\end{propo}

\begin{corol}{}\label{cor1}{\sl Under the conditions of Lemma \ref{lem1},
suppose that the condition (\ref{eq14}) take a place instead of
(\ref{eq4}) and (\ref{eq11}), and the condition
\begin{equation}\label{eq16}
\lim\limits_{x\rightarrow
b}\frac{|f(x)|}{\left(\log\frac{1}{|x-b|}\right)^{\beta_n}}=0\,,
\end{equation}
take a place instead of (\ref{eq3}), where
$\beta_n=\left(\frac{\omega_{n-1}}{A}\right)^{1/(n-1)}.$ Then a
point $x=b$ is a removable singularity of $f.$ }
\end{corol}

{\it Proof.\,} Without loss of generality, we can consider that
$b=0.$ By Lemma \ref{lem1} a point $b$ can not to be an essential
singularity of $f.$ Suppose that $b=0$ is a pole of $f.$ Consider
the composition of the mappings $h=g\circ f,$ where
$g(x)=\frac{x}{|x|^2}$ is inversion under the sphere ${\Bbb
S}^{n-1}.$ Note that a mapping $h$ to be a mapping with finite
length distortion, $K_I(x, f)=K_I(x, h)$ and $h(0)=0.$ Moreover, $h$
is bounded in some neighborhood of zero. Thus there exist
$\varepsilon_1>0$ and $R>0$ such that $|h(x)|\le R$ at every
$|x|<\varepsilon_1.$ Now we apply a Proposition \ref{pr3}. By
(\ref{eq15}), $$|h(x)|=\frac{1}{|f(x)|}\le
\frac{\alpha_n(1+R^2)}{\delta}\left\{\frac{\log\frac{1}{\varepsilon_0}}
{\log\frac{1}{|x|}}\right\}^{\beta_n}\,.$$ Consequently,
$$\frac{|f(x)|}{\left\{\log\frac{1}{|x|}\right\}^{\beta_n}}\ge
\frac{\delta}{\alpha_n(1+R^2){\left\{\log\frac{1}{\varepsilon_0}\right\}}^{\beta_n}}\,.$$
However, the last relation contradicts to the (\ref{eq16}). The
contradiction obtained above prove that $b=0$ is a removable
singularity of $f.$ $\Box$

\setcounter{equation}{0}
\section{The proof of the main results} {\bf Proof of the
statement $1^{\,\prime}$} follows from (\ref{eq14}) which holds for
every function $Q\in FMO(b),$ see, for instance, Corollary 2.3 in
\cite{IR}, or Lemma 6.1 of Ch. VI in \cite{MRSY$_2$}, and from the
Lemma \ref{lem1}. $\Box$

\begin{corol}{}\label{cor2} {\sl\, Let
$f:D\setminus\{b\}\rightarrow {\Bbb R}^n$ be an open and a discrete
mapping with finite length distortion. Suppose that there exists a
measurable function such that $Q:D\rightarrow[1, \infty],$ such that
$K_I(x,f)\le Q(x)$ at a.e. $x\in D$ and $Q(x)\in FMO(b).$

There exists $p_0>0$ such that
\begin{equation}\label{eq18} \lim\limits_{x\rightarrow
b}\frac{|f(x)|}{\left(\log\frac{1}{|x-b|}\right)^{p_0}}=0
\end{equation} implies that $b=0$ is a removable singularity of
$f.$}
\end{corol}

{\bf Proof} follows directly from the Statement $1^{\,\prime}$ and
Corollary \ref{cor1}.$\Box$

\medskip In what follows
$q_{x_0}(r)$ denotes the integral average of $Q(x)$ under the sphere
$|x-x_0|=r,$
\begin{equation}\label{eq17}
q_{x_0}(r):=\frac{1}{\omega_{n-1}r^{n-1}}\int\limits_{|x-x_0|=r}Q(x)\,dS\,,
\end{equation}
where $dS$ is element of the square of the surface $S.$

\medskip
\begin{theo}{}\label{th1}{\sl Let $b\in D$ and
$f:D\setminus\{b\}\rightarrow {\Bbb R}^n$ be an open and a discrete
mapping with finite length distortion. Suppose that there exists
$\delta>0$ such that the relation (\ref{eq2}) holds at every $x\in
B(b, \delta)$ and some constants $p>0$ and $C>0.$ Let there exists a
function $Q:D\rightarrow[1, \infty],$ such that $K_I(x,f)\le Q(x)$
a.e. $x\in D$ and $q_{b}(r)\le C\cdot
\left(\log\frac{1}{r}\right)^{n-1}$ as $r\rightarrow 0.$ The a point
$b$ is a pole or a removable singularity of $f.$

Moreover, in addition, if the relation (\ref{eq18}) holds as
$p_0=\left(\frac{1}{C}\right)^{1/(n-1)},$ then a point $b=0$ is
removable for $f.$}\end{theo}

{\it Proof.\,} We can consider that $b=0.$ Let
$\varepsilon_0<\min\left\{\,\,{\rm dist\,}\left(0,
\partial D\right),\quad 1\right\}.$ Set
$\psi(t)=\frac{1}{t\,\log{\frac{1}{t}}}.$ Note that
$$\int\limits_{\varepsilon<|x|<\varepsilon_0}
\frac{Q(x)dm(x)}{\left(|
x|\log{\frac{1}{|x|}}\right)^n}=\int\limits_{\varepsilon}^{\varepsilon_0}
\left(\int\limits_{|x|\,=\,r}\frac{Q(x)dm(x)}{\left(|
x|\log{\frac{1}{|x|}}\right)^n}\, dS\,\right)\,dr\le C\cdot\omega_
{n-1}\cdot I(\varepsilon, \varepsilon_0)\,,$$
where as above $I(\varepsilon,
\varepsilon_0):=\int\limits_{\varepsilon}^{\varepsilon_0}\psi(t) dt
= \log{\frac{\log{\frac{1}
{\varepsilon}}}{\log{\frac{1}{\varepsilon_0}}}}.$ Thus the
conditions (\ref{eq4}) and (\ref{eq11}) of Lemma \ref{lem1} hold at
$\psi$ which was given above. The second statement of the Theorem
\ref{th1} follows from the Corollary \ref{cor1}. $\Box$

\section{Corollaries. The precision of the conditions}
\setcounter{equation}{0} Recall that a point $y_0\in D$ is said to
be a {\it branch point} of the mapping $f:D\rightarrow {\Bbb R}^n,$
if for every neighborhood $U$ of the point $y_0$ a restriction
$f|_{U}$ fails to be a homeomorphism. A set of all branch sets of
$f$ is denoted by $B_f.$ The following statement can be found as
Theorem 1 in \cite{Sev$_2$}. Let $f:D\rightarrow {\Bbb R}^n$ be an
open discrete mapping of the class $W_{loc}^{1,n}(D)$ such that
either $K_O(x,f)\in L_{loc}^{n-1},$ or $K_I(x,f)\in L_{loc}^1,$ and
$m(B_f)=0.$ Then $f$ to be a mapping with finite length distortion.
Founded on the Statement $1^{\,\prime}$ and on Theorem \ref{th1}, we
have the following.

\begin{theo}{}\label{th2}{\sl Let $b\in D$ and $f:D\setminus\{b\}\rightarrow {\Bbb R}^n$ be an open
discrete mapping of the class $W_{loc}^{1,n}(D),$ for which either
$K_O(x,f)\in L_{loc}^{n-1},$ or $K_I(x,f)\in L_{loc}^1,$ and
$m(B_f)=0.$ Suppose that there exists $\delta>0$ such that the
inequality
$$|f(x)|\le C \left(\log\frac{1}{|x-b|}\right)^{p}$$
take a place for all $x\in B(b, \delta)$ and some constants $p>0$
and $C>0.$ Besides that, suppose that there exists a measurable
function $Q:D\rightarrow[1, \infty]$ such that $K_I(x,f)\le Q(x)$ at
a.e. $x\in D$ and $Q(x)\in FMO(b).$ Then a point $b$ is a removable
singularity, or a pole of $f.$

Moreover, there exists a number $p_0>0$ such that the condition
$$\lim\limits_{x\rightarrow
b}\frac{|f(x)|}{\left(\log\frac{1}{|x-b|}\right)^{p_0}}=0$$ implies
that a point $b$ to be a removable singularity of $f.$ }\end{theo}

\begin{theo}{}\label{th3}{\sl\, All of the conclusions of the Theorem \ref{th2} take a place if
the assumption $Q(x)\in FMO(b)$ to replace on the requirement
$q_{b}(r)\le C\cdot \left(\log\frac{1}{r}\right)^{n-1}$ as
$r\rightarrow 0.$ In this case, we can take
$p_0=\left(\frac{1}{C}\right)^{1/(n-1)}.$} \end{theo}

\medskip
The following result shows that the conditions on $Q$ which are
given above can not to be done more weaker, for instance, we can not
replace it by the assumption $Q\in L_{loc}^q,$ $q\ge 1,$ for every
sufficiently large $q.$

\medskip
\begin{theo}{}\label{th4}{\sl\, Given $p>0$ and $q\in [1, \infty),$ there exists a
homeomorphism $f:{\Bbb B}^n\setminus\{0\}\rightarrow {\Bbb R}^n$
with finite length distortion which belongs to $f\in W_{loc}^{1,n},$
and $f^{-1}\in W_{loc}^{1,n},$ such that $K_I\in L^q_{loc}({\Bbb
B}^n),$
\begin{equation}\label{eq19}|f(x)|\le 2
\left(\log\frac{1}{|x|}\right)^{p}
\end{equation}
at every $x\in B(0, 1/e)\setminus\{0\},$ and a point $b=0$ to be an
essential singularity of $f.$ Moreover, $f$ is a bounded mapping in
this case.}\end{theo}

{\it Proof.\,} The desired homeomorphism $f: {\Bbb
B}^n\setminus\{0\}\rightarrow {\Bbb R}^n$ can be given as
$$f(x)=\frac{1+|x|^{\alpha}}{|x|}\cdot x\,,$$
where $\alpha\in \left(0, n/q(n-1)\right).$ We can consider that
$\alpha<1.$ Note that $f$ maps ${\Bbb B}^n\setminus\{0\}$ onto the
ring $\{1<|y|<2\}$ in ${\Bbb R}^n,$ and the cluster set
$C(f,0)=\{|y|=1\}.$ In particular, it follows from here that $x_0=0$
is an essential singularity of $f.$ It is clear that $f\in C^1({\Bbb
B}^n\setminus\{0\})$ and, consequently, $f\in W_{loc}^{1,n},$
moreover, $K_I(x,f)=\left(\frac{1+|x|^{\,\alpha}}{\alpha
|x|^{\,\alpha}}\right)^{\,n-1}\le\frac{C}{|x|^{\,(n-1)\alpha}},$ see
Proposition 6.3 of Ch. VI in \cite{MRSY$_2$}. Thus, $K_I(x,f)\in
L^q({\Bbb B}^n)$ because $\alpha (n-1) q<n.$ Besides that, to note
that $f$ is locally quasiconformal mapping and hence $f^{\,-1}\in
W_{loc}^{1,n}.$ Thus, $f$ is a mapping with finite length distortion
in ${\Bbb B}^n\setminus\{0\}$ by Theorem 4.6 in \cite{MRSY$_1$}, see
also Theorem 8.1 of Ch. VIII in \cite{MRSY$_2$}.

The mapping $f$ is bounded in ${\Bbb B}^n\setminus\{0\},$ in
particular, $f$ satisfies the inequality $|f(x)|\le 2$ at $x\in
{\Bbb B}^n\setminus\{0\}.$ From other hand, a function
$s(x):=\left(\log\frac{1}{|x|}\right)^{q}$ satisfies $|s(x)|\ge 1$
for all $|x|\le 1/e.$ From here we have a relation (\ref{eq19}).

Thus, we construct a mapping $f$ which have an essential isolated
singularity, and satisfying all of the conditions of the Theorem
\ref{th4}. $\Box$

\medskip
The following statement shows that the condition of the openness of
the mapping $f$ is essential.

\medskip
\begin{theo}{}\label{th5}{\sl\, There exist a discrete mapping
$f:{\Bbb R}^n\setminus\{0\}\rightarrow {\Bbb R}^n$ with finite
length distortion such that $K_I\equiv 1,$ satisfying to the
condition
\begin{equation}\label{eq20}|f(x)|\le \left(\log\frac{1}{|x|}\right)^{p}
\end{equation}
at every $x\in B(0, 1/e)\setminus\{0\}$ and $p>0,$ such that a point
$b=0$ is an essential singularity of $f.$ }\end{theo}

{\it Proof.\,} Consider the division of ${\Bbb R}^n$ by the cubes
$$C_{k_1,\ldots,k_n}=\prod\limits_{i=1}^{n}\left[2k_i-1, 2k_i+1\right]\,,
\quad k_i\in {\Bbb Z}\,.$$
Consider a cube $C_{k_1,\ldots,k_n}$ with $k_1,\ldots,k_n \ge 0;$
the case of the different signs of $k_i$ can be considered by
analogy. Let $x=(x_1,\ldots, x_n)\in C_{k_1,\ldots,k_n}.$ If
$k_1=0,$ $g_{m_1}:={\rm id}.$ Let $k_1>0.$ Set
$f_{1,\ldots,1,1}(x)=y_{1,\ldots,1},$ where $y_{1,\ldots,1,1}$ be a
symmetric reflection of $x$ under the hyperplane $x_1=2k_1-1.$ If
$2k_1-3=-1,$ the process is finished. Let $2k_1-3>-1,$ then
$f_{1,\ldots,1,2}(x)=y_{1,\ldots,1,2},$ where $y_{1,\ldots,1,2}$ be
a symmetric reflection of the point $y_{1,\ldots,1}$ under the
hyperplane $x_1=2k_1-3.$ If $2k_1-5=-1,$ the process is finished. In
other case we continue, $f_{1,\ldots,1,3}(x)=y_{1,\ldots,1,3}.$ Etc.
After a finite number of the steps $m_1$ we have a mapping
$g_{m_1}=f_{1,\ldots,1,m_1}\circ\cdots\circ f_{1,\ldots,1,1},$ such
that $g_{m_1}(x)\in C_{0,k_2,k_3\ldots,k_n}.$

Follow, if $k_2=0,$ then $g_{m_2}:=g_{m_1}.$ As $k_2>0,$ we repeat
the above transformations with the coordinate $x_2$ and the point
$x_{m_1}:=g_{m_1}(x).$ Set $f_{1,\ldots,1,2, m_1}(x)$
$=y_{1,\ldots,1,2,m_1},$ where $y_{1,\ldots,1,2,m_1}$ be a symmetric
reflection of the point $x_{m_1}$ under the hyperplane $x_2=2k_2-1.$
If $2k_2-3=-1,$ the process is finished. In other case we continue.
Now, we have a mapping $g_{m_2}=f_{1,\ldots,m_2,m_1}\circ\cdots\circ
f_{1,\ldots,2,m_1},$ such that $g_{m_2}(x_{m_1})\in
C_{0,0,k_3\ldots,k_n}.$

Etc. After some number of the steps $m_0=m_1+m_2+\ldots+m_n$ we
obtain a mapping $G_0=g_{m_n}\circ g_{m_{n-1}}\circ\cdots
g_{m_{2}}\circ g_{m_{1}},$ such that the image $x_{m_n}$ of the
point $x$ under the mapping $G_0$ lies in the cube
$C_{0,0,0\ldots,0}.$ The compressing $G_1(x)=\frac{\sqrt{n}}{n}\cdot
x$ maps $C_{0,0,0\ldots,0}$ into some cube $A_0,$ which lies in
$\overline{{\Bbb B}^n}.$ Set $G_2:=G_1\circ G_0.$

Note that a point $z_0=\infty$ is an essential singularity of $G_2,$
moreover, $C(G_2,\infty)=A_0\subset \overline{{\Bbb B}^n}.$ Then a
mapping
\begin{equation}\label{eq25*}
g:=G_2\circ G_3\,,
\end{equation}
$G_3(x)=\frac{x}{|x|^2},$ has an essential singularity $b=0,$ and
%
%
$$C(g, 0)\subset \overline{{\Bbb B}^n}\,.$$
%
%
By the construction of $G_2,$ which is given by (\ref{eq25*}), $G_2$
is a discrete mapping, preserves the lengths of curves in ${\Bbb
R}^n,$ is differentiable a.e. and has $N$ and $N^{-1}$ --
properties. Thus, $g$ is a mapping with finite length distortion,
moreover, it is easy to see that $K_I(x,g)=1.$ Finally, $|g(x)|\le
1$ at every $x\in {\Bbb R}^n\setminus\{0\}.$ Thus, (\ref{eq20})
holds at every $x\in B(0, 1/e)\setminus\{0\}.$

The desired mapping have been constructed. $\Box$


\large
{\bf \noindent Evgenii A. Sevost'yanov} \\
Institute of Applied Mathematics and Mechanics,\\
National Academy of Sciences of Ukraine, \\
74 Roze Luxemburg str., 83114 Donetsk, UKRAINE \\
Phone: +38 -- (062) -- 3110145, \\
Email: brusin2006@rambler.ru

\begin{thebibliography}{99}

\bibitem[AIM]{AIM} K.~Astala, T.~Iwaniec and G.~Martin, {\it Elliptic
Partial Differential Equations and Quasiconformal Mappings in the
Plane}, Princeton: Princeton University Press, 2009.

\bibitem[BGMV]{BGMV} C.J.~Bishop, V.Ya.~Gutlyanskii, O.~Martio, M.~Vuorinen, $"$On conformal dilatation in
space$"$, {\it Intern. J. Math. and Math. Scie.}, \textbf{22}
(2003), 1397--1420.

\bibitem[BI]{BI} B.~Bojarski and  T.~Iwaniec, $"$Analytical foundations
of the theory of quasiconformal mappings in ${\Bbb R}^n$$"$, {\it
Ann. Acad. Sci. Fenn. Ser. A 1 Math.}, \textbf{8}:2 (1983),
257--324.

\bibitem[IR]{IR} A.~Ignat'ev and V.~Ryazanov, $"$Finite mean oscillation
in the mapping theory$"$, {\it Ukrainian Math. Bull.}, \textbf{2}:3
(2005), 403--424.

\bibitem[IM]{IM}  T.~Iwaniec and G.~Martin, {\it Geometrical Function
Theory and Non--Linear Analysis}, Oxford: Claren\-don Press, 2001.

\bibitem[JN]{JN} F.~John, L.~Nirenberg, $"$On functions of bounded
mean oscillation$"$, {\it Comm. Pure Appl. Math.}, \textbf{14}
(1961), 415--426.

\bibitem[MRSY$_1$]{MRSY$_1$} O.~Martio, V.~Ryazanov, U.~Srebro, E.~Yakubov,
$"$Mappings with finite length distortion$"$, {\it J. d'Anal.
Math.}, \textbf{93} (2004), 215--236.

%
\bibitem[MRSY$_2$]{MRSY$_2$} O.~Martio, V.~Ryazanov, U.~Srebro, E.~Yakubov,
{\it Moduli in Modern Mapping Theory}, New York: Springer Science +
Business Media, LLC, 2009.

\bibitem[Mikl]{Mikl} V.M.~Miklyukov, $"$The relation distance of M.A. Lavrent'ev and prime ends on non--parametric
surface$"$, {\it Ukr. Math. Bull.}, \textbf{1}:3 (2004), 353--376.

\bibitem[Pol]{Pol} E.A.~Poletskii, $"$The modulus method for
non--homeomorphic quasiconformal mappings$"$, {\it Math. Sb.},
\textbf{83}:2 (1970), 261--272 (in Russian).

\bibitem[Po]{Po} S.P.~Ponomarev, $"$The $N^1$ -- property of maps and
Luzin’s condition $(N)$$"$, {\it Math. Notes}, \textbf{58}:3 (1995),
960-–965.

\bibitem[Re$_1$]{Re$_1$} Yu.G.~Reshetnyak, $"$Generalized derivatives an differentiability
a.e.$"$, {\it Math. Sb.}, \textbf{75}:3  (1968), 323--334 (in
Russian).

\bibitem[Re$_2$]{Re$_2$} Yu.G.~Reshetnyak, {\it Space mappings with bounded distortion},
English transl.: Amer. Math. Soc., Providence, RI, 1989.

\bibitem[Ri]{Ri} S.~Rickman, {\it Quasiregular mappings},
Berlin etc.: Springer-Verlag, 1993.

\bibitem[Sal]{Sal} R.~Salimov, $"$On regular homeomorphisms in the
plane$"$, {\it Ann. Acad. Sci. Fenn. Ser. A1. Math.,} \textbf{35}
(2010), 285--289.

\bibitem[UV]{UV} A.D.~Ukhlov and S.K.~Vodop'yanov, $"$Sobolev spaces and mappings with bounded $(P; Q)$ --
distortion on Carnot groups$"$, {\it Bull. Sci. Mat.}, \textbf{54}:4
(2009), 349--370.

\bibitem[Va$_1$]{Va$_1$} J.~V\"{a}is\"{a}l\"{a}, {\it Lectures on $n$ --
Dimensional Quasiconformal Mappings}, Lecture Notes in Math.,
\textbf{229}, Berlin etc.: Springer--Verlag, 1971.

\bibitem[Va$_2$]{Va$_2$} J.~V\"{a}is\"{a}l\"{a}, $"$Modulus and capacity
inequalities for quasiregular mappings$"$, {\it Ann. Acad. Sci.
Fenn. Ser. A 1 Math.}, \textbf{509} (1972), 1--14.

\bibitem[Sev$_1$]{Sev$_1$} E.A.~Sevost'yanov, $"$On the normality of space
mappings with branching$"$, {\it Ukr. Math. J.}, \textbf{60}:10
(2008), P. 1618--1632.

\bibitem[Sev$_2$]{Sev$_2$} E.A.~Sevost'yanov, $"$Generalization of one Poletskii lemma to
classes of space mappings$"$, {\it Ukr. Math. J.}, \textbf{61}:7
(2009), P. 1151--1157.

\bibitem[Sev$_3$]{Sev$_3$} E.A.~Sevost'yanov, $"$The V\"{a}is\"{a}l\"{a}
inequality for mappings with finite length distortion$"$, {\it
Complex Variables and Elliptic Equations}, \textbf{55}: 1--3 (2010),
91--101.

\bibitem[Sev$_4$]{Sev$_4$} E.A.~Sevost'yanov,
$"$Towards a theory of removable singularities for maps with
unbounded characteristic of quasi-conformity$"$,  {\it Izvestiya:
Mathematics}, \textbf{74}:1 (2010), 151--165.

\end{thebibliography}
\end{document}